%

\magnification=\magstep1
\def\forces{\parallel\!\!\! -}


\def\hexnumber#1{\ifcase#1 0\or1\or2\or3\or4\or5\or6\or7\or8\or9\or
	A\or B\or C\or D\or E\or F\fi }

\font\teneuf=eufm10
\font\seveneuf=eufm7
\font\fiveeuf=eufm5
\newfam\euffam
\textfont\euffam=\teneuf
\scriptfont\euffam=\seveneuf
\scriptscriptfont\euffam=\fiveeuf


\font\tenmsx=msam10
\font\sevenmsx=msam7
\font\fivemsx=msam5
\font\tenmsy=msbm10
\font\sevenmsy=msbm7
\font\fivemsy=msbm5
\newfam\msxfam
\newfam\msyfam
\textfont\msxfam=\tenmsx  \scriptfont\msxfam=\sevenmsx
  \scriptscriptfont\msxfam=\fivemsx
\textfont\msyfam=\tenmsy  \scriptfont\msyfam=\sevenmsy
  \scriptscriptfont\msyfam=\fivemsy
\edef\msx{\hexnumber\msxfam}

\mathchardef\upharpoonright="0\msx16
\let\restriction=\upharpoonright
\def\Bbb#1{\tenmsy\fam\msyfam#1}

\def\re{{\restriction}}

\def\Smallskip{\vskip1.4truecm}
\def\Bigskip{\vskip2.2truecm}

\def\qed{{\vcenter{\hrule height.4pt \hbox{\vrule width.4pt height5pt
 \kern5pt \vrule width.4pt} \hrule height.4pt}}}
\def\ok{\vbox{\hrule height 8pt width 8pt depth -7.4pt
    \hbox{\vrule width 0.6pt height 7.4pt \kern 7.4pt \vrule width 0.6pt height 7.4pt}
    \hrule height 0.6pt width 8pt}}
\def\nt{{\leq}\kern-1.5pt \vrule height 6.5pt width.8pt depth-0.5pt \kern 1pt}
\def\sd{{\times}\kern-2pt \vrule height 5pt width.6pt depth0pt \kern1pt}
\def\zp#1{{\hochss Y}\kern-3pt$_{#1}$\kern-1pt}

\def\BB{{\Bbb B}}
\def\CC{{\Bbb C}}

\def\B{{\cal B}}

\def\F{{\cal F}}

\def\I{{\cal I}}
\def\J{{\cal J}}

\def\M{{\cal M}}
\def\N{{\cal N}}

\def\Add#1{{\sanse add}$({\cal #1})$}
\def\Cov#1{{\sanse cov}$({\cal #1})$}
\def\Unif#1{{\sanse unif}$({\cal #1})$}
\def\Cof#1{{\sanse cof}$({\cal #1})$}
\def\Coll#1#2{{\sanse Coll}$(#1,#2)$}
\def\Cove#1#2{{\sanse cov}$^\star (#1 , #2)$}

\def\sm{{\smallskip}}
\def\ce#1{{\centerline{#1}}}

\def\no{{\noindent}}
\def\la{{\langle}}
\def\ra{{\rangle}}
\def\sub{\subseteq}

\def\em{{\emptyset}}
\def\sem{\setminus}
\def\omom{{\omega^\omega}}

\def\twoom{{2^\omega}}

\def\Lolear{\Longleftarrow}
\def\Loriar{\Longrightarrow}

\font\small=cmr8 scaled\magstep0
\font\smalli=cmti8 scaled\magstep0
\font\capit=cmcsc10 scaled\magstep0
\font\capitg=cmcsc10 scaled\magstep1

\font\dunhgg=cmdunh10 scaled\magstep2

\font\sanse=cmss10 scaled\magstep0

\font\bolds=cmssdc10 scaled\magstep0

\overfullrule=0pt
\openup1.5\jot

\ce{}
\Smallskip
\ce{\dunhgg  Nicely generated and chaotic ideals}
\footnote{}{{\openup-6pt {\small {\smalli
1991 Mathematics subject classification.}
03E05 03E40 28A05 54H05    \par
{\smalli Key words and phrases.} Meager sets, null sets,
Cohen reals, random reals.
\endgraf}}}
\Bigskip
\ce{\capitg J\"org Brendle\footnote{$^\star$}
{{\small Supported by DFG--grant Nr. Br 1420/1--1.}}}
\Smallskip
\no Mathematisches Institut der Universit\"at
T\"ubingen, Auf der Morgenstelle 10, 72076
T\"ubingen, Germany; email: {\sanse
jobr@michelangelo.mathematik.uni--tuebingen.de}
\Bigskip
\ce{\capit Abstract}
\bigskip
\no We show that if the real line is the disjoint union
of $\kappa$ meager sets such that every meager set is
contained in a countable union of them, then $\kappa =
\omega_1$. This answers a question addressed by Jacek
Cicho\'n. We also prove two theorems saying roughly that
any attempt to produce the isomorphism type of the meager
ideal in the Cohen real and the random real extensions
must fail. All our results hold for {\it meager} replaced
by {\it null}, as well.
\vfill\eject

\ce{\capitg Introduction}
\Smallskip
The purpose of this note is to prove three results concerning
the structure of the $\sigma$--ideals of meager sets $\M$
and of null sets $\N$ on the reals. By {\it the reals} we shall
mean the Cantor space $\twoom$, equipped with the product
topology and the product measure (where $2$ carries the 
discrete topology and the measure giving both $0$ and $1$
measure ${1\over 2}$). \par
In the first section we investigate ideals which are generated
by small unions of sets in a partition of the underlying
set. We show that --- except for trivial cases --- the
meager and null ideals are not of this form (Corollary 6).
In the second section we study the structure of the meager
ideal $\M$ and the null ideal $\N$ in both the Cohen and
random real models. Unfortunately, all our results (Theorems
9, 9$^\star$, 10, 10$^\star$) are on the non--structural
side. \par
{\it Notational remarks.} Given cardinals $\kappa$ and $\lambda$,
$[\lambda]^{<\kappa}$ denotes the set of subsets of $\lambda$
of size $<\kappa$ and $[\lambda]^\kappa$ stands for the set of
subsets of $\lambda$ of size $\kappa$. Let $\B$ denote the
Borel subsets of $\twoom$. Then $\CC = \B / \M$ is the
{\it Cohen algebra}, and $\BB = \B / \N$ is the {\it random
algebra}. More generally, for a cardinal $\kappa$, $\CC_\kappa$
($\BB_\kappa$, respectively) is the algebra adding $\kappa$ Cohen
(random, resp.) reals; and for $A \sub\kappa$, $\CC_A$ ($\BB_A$,
resp.) is the complete subalgebra adjoining the $|A|$ many
Cohen (random, resp.) reals with index in $A$. When talking about
an ideal $\I$ on a set $X$ (i.e. $\I \sub P(X)$), we assume
throughout that $\bigcup \I = X$ and $X \notin \I$ to avoid
pathologies. For more set--theoretical notation we refer
the reader to [Je] or [Ku 1]. For the basic facts about 
Cohen and random forcing see [Ku 2].
\Bigskip

\ce{\capitg 1. Nicely generated ideals}
\Smallskip

{\capit Definition 1.} Let $X$ be a set, and let $\kappa \leq \lambda$
be cardinals with $\lambda \leq |X|$. An ideal $\I \sub P(X)$
is {\it $(\kappa , \lambda)$--nicely generated} if there is a
partition $\la X_\alpha \in \I ; \; \alpha < \lambda \ra$ of $X$
such that \par
\item{(i)} $\bigcup_{\alpha \in\Delta} X_\alpha \in \I$ for
any $\Delta \in [\lambda]^{<\kappa}$ and \par
\item{(ii)} for any $Y \in \I$ there is $\Delta \in [\lambda]^{<\kappa}$
with $Y \sub \bigcup_{\alpha \in \Delta} X_\alpha$. $\qed$ \par
\sm

\no In this section we try to characterize the situations in
which $\M$ and $\N$ can be nicely generated.
\sm

{\capit Definition 2.} For an ideal $\I \sub P(X)$ we define the
following cardinal characteristics.
\sm
\Add I $=\min \{ |\F| ; \; \F \sub \I \;\land\; \bigcup \F \notin \I \}$;
\par
\Cov I $=\min \{ |\F| ; \; \F \sub \I \;\land\; \bigcup \F = X \}$;
\par
\Unif I $=\min \{ |Y| ; \; Y \sub X \;\land\; Y\notin\I \}$;
\par
\Cof I $=\min \{ |\F| ; \; \F \sub\I \;\land\; \forall Y\in\I \;\exists
Z\in\F \; (Y \sub Z) \}$. $\qed$
\sm

\no It is well--known which inequalities are provable in $ZFC$
between these cardinal invariants for $\I = \M$ and $\I = \N$ 
(see [Fr], [BJ]). Furthermore, their values in various
$ZFC$--models have been studied intensively (see [BJS], [BJ]).
\sm

{\it Fact 3. If $X , \kappa , \lambda$ are as in Definition 1, and
$\I \sub P(X)$ is $(\kappa , \lambda)$--nicely generated, then
$$\eqalign{\hbox{\Add I} &= cf(\kappa), \cr
\hbox{\Cov I} &=\cases{\lambda & if $\kappa < \lambda$, \cr
                       cf(\lambda) & if $\kappa = \lambda$, \cr} \cr
\hbox{\Unif I} &= \kappa, \cr
\hbox{\Cof I} &= \hbox{\Cove {<\kappa} \lambda}, \cr}$$
where \Cove {<\kappa} \lambda $= \min \{ |\F| ; \; \F \sub [\lambda]^{
< \kappa} \;\land\; \forall Y \in [\lambda]^{<\kappa} \;
\exists Z \in \F \; (Y \sub Z) \}$.} $\qed$
\sm

\no The following result is well--known. We include a proof
for completeness' sake.
\sm

{\capit Lemma 4.} (Folklore) {\it Assume $\{ A_\alpha ; \; \alpha < \lambda^+
\}$ is a family of sets of size $< \lambda$. Then there are $\Lambda
\in [\lambda^+]^{\lambda^+}$ and a set $A$ of size $\leq\lambda$
such that $A_\alpha \cap A_\beta \sub A$ for distinct
$\alpha , \beta \in \Lambda$. In case $\lambda$ is regular,
we can find such $\Lambda$ and $A$ with $|A| < \lambda$.}
\sm
{\it Proof.} Without loss $A_\alpha \sub \lambda^+$ for $\alpha \in
\lambda^+$. Let $\beta_\alpha < \lambda$ be the order type
of $A_\alpha$ under the inherited ordering. Find $\Gamma \in
[\lambda^+]^{\lambda^+}$ and $\beta < \lambda$ with $\beta_\alpha
= \beta$ for $\alpha \in \Gamma$. Let $\gamma < \beta$ be minimal
such that $| \{ \delta ; \; \exists \alpha \in \Gamma \; ( \delta$
is the $\gamma$--th element of $A_\alpha) \} | = \lambda^+$ if there
is such a $\gamma$; otherwise put $\gamma = \beta$. If $\gamma < 
\beta$ let $\delta_\alpha$ be the $\gamma$--th element of $A_\alpha$
for $\alpha\in\Gamma$, and put $B_\alpha = A_\alpha \sem \delta_\alpha$;
otherwise $\delta_\alpha = \lambda^+$. It
is now easy to find $\Delta \in [\Gamma]^{\lambda^+}$ such that
$B_\alpha \cap B_\beta = \em$ for distinct $\alpha , \beta \in
\Delta$. Set $C = \bigcup_{\alpha \in \Delta} A_\alpha \cap \delta_\alpha$.
In the general case, $\Lambda = \Delta$ and $A=C$ satisfy the
requirements of the Lemma. In case $\lambda$ is regular and
$|C|  =\lambda$, we easily find $A \sub C$ with $|A|
< \lambda$ and $\Lambda \in [\Delta]^{\lambda^+}$ with $A_\alpha
\cap \delta_\alpha \sub A$ for $\alpha \in \Lambda$. $\qed$
\sm

\no Notice that in case $\lambda^{<\lambda} = \lambda$, the
$\Delta$--system lemma says something much stronger [Ku 1, II.1.6]. ---
We are ready to state and prove the main result of this section. 
\sm

{\bolds Theorem 5.} {\it Let $\omega_1 \leq\kappa\leq\lambda\leq\twoom$
be cardinals, and assume the meager ideal $\M$ on $\twoom$
is $(\kappa , \lambda)$--nicely generated. Then $\kappa = \lambda$.}
\sm
{\it Proof.} We assume $\kappa < \lambda$ and we seek a contradiction.
Without loss $\lambda = \kappa^+ = \twoom$. Otherwise we force
with the collapse \Coll {\kappa^+} \twoom $= \{ p ; \; |p|
\leq \kappa$ and $p: \kappa^+ \to \twoom$ is a partial function$\}$.
The extension satisfies $|\lambda| = |\twoom| = \kappa^+$,
and $\M$ is $(\kappa , \kappa^+)$--nicely generated (because
no reals and no Borel sets are adjoined).
\par
Let $\la X_\alpha ; \; \alpha < \kappa^+ \ra \sub \twoom$ witness
that $\M$ is $(\kappa , \kappa^+)$--nicely generated. Let
$\la M_\alpha ; \; \alpha < \kappa^+ \ra$ be sequence of models
of $ZFC$ satisfying $|M_\alpha| = \kappa$, $M_\alpha \sub 
M_\beta$ for $\alpha \leq \beta$, and $\bigcup_{\alpha < \kappa^+}
M_\alpha \supseteq \twoom$, such that there are reals $c_\alpha
\in M_{\alpha+1}$ ($\alpha \in\kappa^+$) Cohen over $M_\alpha$.
Such a sequence exists as \Cov M $=\kappa^+$ (by Fact 3).
\par
For a real $x\in\twoom$, let $A_x = \{ f\in\twoom ; \;
\forall n \; (x(2n) = f(2n)) \}$ and $B_x = \{ f\in\twoom ; \;
\forall n \; (x(2n+1) = f(2n+1)) \}$. The $A_x$'s and $B_x$'s
are closed nowhere dense sets. Recall that if $M \sub N$ are
models of $ZFC$, $c\in N$ is $\CC$--generic over $M$ and $d$
is $\CC$--generic over $N$, then $(c,d)$ is $\CC\times\CC$--generic
over $M$, and thus $f = f(c,d)$ defined by $f(2n) = c(2n)$
and $f(2n + 1) = d(2n + 1)$ is $\CC$--generic over $M$ as well (see
[Ku 2] for details).
\par
We put $A_\alpha = A_{c_\alpha}$, $B_\alpha = B_{c_\alpha}$ and
$D_\alpha = \{ \beta < \kappa^+ ; \; A_\alpha \cap X_\beta \neq\em \}$
for $\alpha < \kappa^+$. By assumption $|D_\alpha| < \kappa$, and
we can apply Lemma 4 to the family $\{ D_\alpha ; \; \alpha \in
\kappa^+ \}$ to get $\Lambda \in [\kappa^+]^{\kappa^+}$ and $D$ of
size $\leq \kappa$ with $D_\alpha \cap D_\beta \sub D$ for distinct
$\alpha , \beta \in \Lambda$. For $\alpha \in D$ find a Borel
meager set $S_\alpha$ containing $X_\alpha$. The $S_\alpha$
are coded by reals, so there is $\beta_0 \in \kappa^+$
such that $M_{\beta_0}$ contains all $S_\alpha$ where $\alpha \in D$. Let
$\la \beta_\alpha ; \; 0 < \alpha \leq \kappa \ra$ be a
strictly increasing sequence of elements of $\Lambda$ larger than
$\beta_0$. Let $f_\alpha = f(c_{\beta_\alpha} , c_{\beta_\kappa} )$
for $0 \leq \alpha < \kappa$. Then $\{ f_\alpha ; \;
\alpha < \kappa \} \sub B_{\beta_\kappa}$. Also $f_\alpha \in
A_{\beta_\alpha}$, and thus $f_\alpha \in X_{\gamma_\alpha}$
for some $\gamma_\alpha \in D_{\beta_\alpha}$. As $f_\alpha$ is Cohen over
$M_{\beta_0}$, $f_\alpha \notin S_\gamma$ --- and in particular $f_\alpha
\notin X_\gamma$ --- for $\gamma \in D$. Hence $\gamma_\alpha \in
D_{\beta_\alpha} \sem D$, and therefore the $\gamma_\alpha$'s
are distinct for distinct $\alpha$'s. This entails that
$B_{\beta_\kappa}$ is not included in the union of less than
$\kappa$ many $X_\delta$'s, a contradiction. $\qed$
\sm

{\bolds Theorem 5$^\star$.} {\it Let $\omega_1 \leq \kappa
\leq \lambda \leq \twoom$ be cardinals, and assume the null ideal
$\N$ on $\twoom$ is $(\kappa , \lambda)$--nicely generated.
Then $\kappa = \lambda$.}
\sm
{\it Proof.} Simply replace all instances of {\it Cohen}, $\CC$,
and {\it meager} by {\it random}, $\BB$, and {\it null} in the
above proof --- and notice that the sets $A_x$ and $B_x$ defined
there are null as well. $\qed$
\sm

{\bolds Corollary 6.} {\it Assume $\kappa$ is regular. Then the
meager ideal $\M$ on $\twoom$ is $(\kappa , \lambda)$--nicely generated 
for some $\lambda \geq \kappa$ iff \Add M $=$ \Cof M $=\kappa$.
A dual statement holds for the null ideal $\N$.}
\sm
{\it Proof.} $(\Loriar)$. Theorem 5 (5$^\star$) and Fact 3. \par
$(\Lolear)$. Any ideal $\I$ with \Add I $=$ \Cof I $= \kappa$
is easily seen to be $(\kappa , \kappa)$--nicely generated. $\qed$
\sm

\no We do not know what happens in the singular case. It might
be consistent that for some singular $\kappa$, $\M$ is $(\kappa ,
\kappa)$--nicely generated in which case we would have \Add M $=$
\Cov M $=cf(\kappa)$, \Unif M $=\kappa$, and \Cof M $=$ \Cove {<\kappa}
\kappa $>\kappa$ (Fact 3). We conjecture, however, that this is not
the case.
\par
Our result answers a question addressed by Jacek Cicho\'n
(see [Mi, problem 15.2]). In case $\kappa = \omega_1$
Theorem 5 (or Corollary 6) shows that the real line cannot be the
disjoint union of $\omega_2$ meager sets such that every meager
set is contained in a countable union of them.
\Bigskip

\ce{\capitg 2. Chaotic ideals}
\Smallskip

{\capit Definition 7.} Let $X$, $Y$ be sets and $\I \sub P(X)$,
$\J \sub P(Y)$ be ideals. We say $\I$ and $\J$ are {\it isomorphic}
($\I \cong \J$) if there is a bijection $f : X \to Y$ such that
$A\in\I$ iff $f(A) \in \J$ for all $A\sub X$. $\qed$
\sm

\no Let $\I_0 $ be the ideal on $\omega_2 \times \omega_2$
generated by rectangles of the form $A \times \omega_2$ where
$A \sub \omega_2$ is countable. This is a typical example
of an $(\omega_1 , \omega_2)$--nicely generated ideal.
Originally, Cicho\'n conjectured that adding $\omega_2$
Cohen reals to a model of $CH$ forces $\M$ to be isomorphic
to $\I_0$. This is false by Theorem 5. However one may still
ask whether there is a {\it reasonably nice} $\sigma$--ideal
on $\omega_2$ (or $\omega_2 \times \omega_2$ or some other underlying
set of size $\omega_2$) such that $\M$ is isomorphic to $\I$
in the Cohen model. Of course, what is meant by {\it reasonably
nice} should be clarified. One way of doing this goes by requiring
that such an ideal exist (i.e. be definable) in any
$ZFC$--model (as does the above example $\I_0$). Certainly,
it won't be the same in different models; however, when stepping
from some model $V$ into a $ccc$ forcing extension $W$, any
$\sigma$--ideal of $V$ still generates a $\sigma$--ideal in $W$.
This leads to the following notion.
\sm

{\capit Definition 8.} Let $V \sub W$ be models of $ZFC$ such that
$W$ is a $ccc$ forcing extension of $V$. Let $X \in W$, and let
$\J \sub P(X)$ with $\J \in W$ be a $\sigma$--ideal.
We say $\J$ is {\it chaotic over $V$} if for any $\sigma$--ideal
$\I \sub P(|X|^W)$ with $\I \in V$ we have $\J \not\cong \I^W$,
where $\I^W$ denotes the ideal generated by $\I$ in $W$. $\qed$
\sm

{\bolds Theorem 9.} {\it Let $\kappa \geq \omega_2$. Then
\sm
\ce{$\forces_{\CC_\kappa} `` \M$ is chaotic over $V "$.}}
\sm
{\it Proof.} Let $\I$ be a $\sigma$--ideal on $\lambda$ in $V$ where
$\forces_{\CC_\kappa} `` \twoom = \lambda "$. We have to show that the
ideal generated by $\I$ in the generic extension is not isomorphic to
$\M$. Assume to the contrary, and let $\dot f$ be a $\CC_\kappa$--name
for a bijection between $\lambda$ and $\twoom$ giving rise
to an isomorphism.
\par
Let $\la Y_n ; \; n\in\omega\ra$ be a partition of $\omega$ into
infinitely many infinite sets. For a real $x\in\twoom$, let
$A_x = \{ f\in \twoom ; \; \exists n\in\omega \; (f \re Y_n =
x \re Y_n ) \}$. The $A_x$'s are meager (and  null). Let
$\dot c_\alpha$ be the name for the $\alpha$--th Cohen real
($\alpha \in \kappa$), and $\dot A_\alpha$ denotes the name for
$A_{\dot c_\alpha}$. For each $\alpha$ we can find a name 
for a meager Borel set $\dot M_\alpha$ and a set $X_\alpha \in \I$
such that
\sm 
\ce{$\forces_{\CC_\kappa} `` \dot A_\alpha \sub \dot M_\alpha \;
\land\; \dot f (X_\alpha) = \dot M_\alpha "$.}
\sm
\no To see this fix $\alpha \in \kappa$. As $\forces_{\CC_\kappa}
`` \dot A_\alpha$ is meager ", find $X_\alpha^0 \in \I$ with
$\forces_{\CC_\kappa} `` \dot A_\alpha \sub \dot f (X_\alpha^0) "$,
using the $ccc$--ness of $\CC_\kappa$ and the fact that $\I$ is
a $\sigma$--ideal. Let $\dot M_\alpha^0$ be a name for a meager Borel
set such that $\forces_{\CC_\kappa} `` \dot f (X_\alpha^0 ) \sub
\dot M_\alpha^0 "$. Iterating this construction find
$X_\alpha^n \in\I$ and names for meager Borel sets $\dot M_\alpha^n$
with $\forces_{\CC_\kappa} `` \dot f (X_\alpha^n ) \sub \dot M_\alpha^n
\sub \dot f (X_\alpha^{n+1}) "$. In the end, put $X_\alpha =
\bigcup_{n\in\omega} X_\alpha^n$ and let $\dot M_\alpha$ be a name
for the union of the $\dot M_\alpha^n$.
\par
As each $\dot M_\alpha$ is essentially a real, we can find
countable sets $\{ \alpha \} \sub D_\alpha \sub \kappa$ such that
the interpretation of $\dot M_\alpha$ lies already in the intermediate
extension via the complete subalgebra $\CC_{D_\alpha}$. Applying
Lemma 4 to $\omega_1$ and the $D_\alpha$, we find $\Lambda \in [\kappa]^{
\omega_2}$ and $D \sub\kappa$ countable ($|D| \leq\omega_1$ suffices)
such that $D_\alpha \cap D_\beta \sub D$ for distinct $\alpha , \beta
\in \Lambda$. Without loss $D \neq\em$.
\par
Step into the intermediate extension $V '$ via $\CC_D$;
and let $\{ \alpha_n ; \; n\in\omega \} \sub \Lambda \sem D$ (in $V$).
For $\beta \in \lambda$, let $x_\beta = \{ n\in\omega ; \; \beta
\in X_{\alpha_n} \}$; $x_\beta$ is (essentially) a real
lying in the ground model $V$. Let $x \in V' \sem V$  
be a subset of $\omega$, and let $\phi\in\omom$ be the
increasing enumeration of $x$. Next let $\dot c$ be the $\CC_{\kappa \sem D}
$--name for the following real:
$$\forces_{\CC_{\kappa \sem D}} `` \dot c \re Y_n = 
\dot c_{\alpha_{\phi(n)}} \re Y_n " .$$
Thus $\dot c$ is forced to be an element of $\dot A_{\alpha_n}$
(and hence of $\dot M_{\alpha_n}$) for $n\in x$. On the other hand,
the interpretation of $\dot c$ is still Cohen over the extension
via $\CC_{D_{\alpha_n} \sem D}$ for $n\notin x$; in particular
$\dot c$ is forced to avoid $\dot M_{\alpha_n}$ for $n\notin x$.
Hence
$$\forces_{\CC_{\kappa \sem D}} `` \dot f^{-1} (\dot c) \in
\bigcap_{n\in x} X_{\alpha_n} \cap \bigcap_{n\notin x} (\lambda
\sem X_{\alpha_n} ) ".$$
This is a contradiction because $\beta \in \bigcap_{n\in x_\beta}
X_{\alpha_n} \cap \bigcap_{n\notin x_\beta} (\lambda \sem X_{\alpha_n} )$
and $x_\beta \neq x$ for all $\beta \in \lambda$. $\qed$
\sm

\no Notice that Theorem 9 provides an alternative argument for
Cicho\'n's original conjecture to be false in the Cohen real model.
We also note that Theorem 9 has the following funny consequence:
assume we add first $\omega_2$ Cohen reals and then again $\omega_2$
Cohen reals; let $\M'$ be the meager ideal in the intermediate
model, and $\M$ the meager ideal in the final model; let $\M^\star$
be the ideal generated by $\M'$ in the second extension. Then
$\M^\star \not\cong \M$, although they have the same cardinal characteristics
and are intuitively similar because they are gotten by the same kind
of extension. --- As in the case of Theorem 5, it is immediate
that the dual statement about the null ideal can be proved in 
exactly the same way.
\sm

{\bolds Theorem 9$^\star$.} {\it Let $\kappa \geq \omega_2$. Then
\sm
\ce{$\forces_{\BB_\kappa} `` \N$ is chaotic over $V "$. $\qed$} }
\sm

\no We now switch to the investigation of the meager ideal in the
random extension (and, dually, the investigation of the null
ideal in the Cohen extension).
\sm

{\bolds Theorem 10.} {\it Let $\kappa \geq \omega_2$. Then
\sm
\ce{$\forces_{\BB_\kappa} `` \M$ is chaotic over $V "$. }}
\sm
{\it Proof.} We start  as in the proof of Theorem 9 with a
$\sigma$--ideal $\I$ on $\lambda$ in $V$ where $\forces_{\BB_\kappa}
`` \twoom = \lambda "$, and a $\BB_\kappa$--name $\dot f$ for a bijection
between $\lambda$ and $\twoom$.
\par
Next recall that $\BB$ adjoins a meager Borel set containing all
ground model reals [Ku 2]. Let $\dot A_\alpha$ be a $\BB_{\{\alpha\}}
$--name for such a set. We find as before $\BB_\kappa$--names
$\dot M_\alpha$ for meager Borel sets and $X_\alpha \in \I$ such that
\sm
\ce{$\forces_{\BB_\kappa} `` \dot A_\alpha \sub \dot M_\alpha \;
\land\; \dot f (X_\alpha) = \dot M_\alpha "$.}
\sm
\no The next paragraph in the proof of Theorem 9 can be taken over
with $\CC$ replaced by $\BB$: we get the $D_\alpha$'s, $\Lambda
\in [\kappa]^{\omega_2}$ and $D$.
\par
Then we step again into the intermediate extension $V'$ via
$\BB_D$, and choose $\{ \alpha_n ; \; n\in\omega \} \sub \Lambda
\sem D$. The $x_\beta$, $\beta < \lambda$, and $x$ are 
as before. Next we step into the extension $V''$ of $V'$ via
$\BB_{\bigcup_{n\notin x} D_{\alpha_n} \sem D}$. In this
model find a real $d$ which is not in the interpretation of any
$\dot M_{\alpha_n}$, $n\notin x$. On the other hand, $d$ is forced
to lie in $\dot A_{\alpha_n}$ for $n\in x$ (in fact, this is true
for any real of $V''$); in particular it will be in $\dot M_{\alpha_n}$
for $n\in x$. Thus we have again
$$\forces_{\BB_{\kappa \sem D}} `` \dot f^{-1} (\dot d)
\in \bigcap_{n\in x} X_{\alpha_n} \cap \bigcap_{n \notin x}
(\lambda \sem X_{\alpha_n} ) " ,$$
where $\dot d$ is the $\BB_{\kappa \sem D}$--name for $d$.
This is a contradiction as in the proof of Theorem 9. $\qed$
\sm

\no An exactly similar argument yields:
\sm

{\bolds Theorem 10$^\star$.} {\it Let $\kappa \geq \omega_2$.
Then
\sm
\ce{$\forces_{\CC_\kappa} `` \N$ is chaotic over $V "$. $\qed$ }}
\sm

\no We have seen several situations in which $\M$ and $\N$ are
chaotic, but the problem remains whether they can be non--chaotic
in a non--trivial way. Of course, if the $ccc$--extension
forces $MA$ to hold, then $\M$ and $\N$ are not chaotic over
the ground model. However, in this case, \Cof M = \Cof N
= \Add M = \Add N $=\twoom$, and the structure of both ideals is
well--known and trivial. One may still hope that
in other models of $ZFC$ (e.g. the ones gotten from adding
iteratively $\omega_2$ Laver or $\omega_2$ Miller reals
with countable support over a model for $CH$), $\M$ and
$\N$ have a nicer structure. Note in this context
that $ccc$ can be replaced by {\it proper and cardinal--preserving}
in Definition 8.

\Bigskip

\ce{\capitg References}
\Smallskip

\itemitem{[BJ]} {\capit T. Bartoszy\'nski and H. Judah,}
{\it Measure and category: the asymmetry,} forthcoming book.
\smallskip
\itemitem{[BJS]} {\capit T. Bartoszy\'nski, H. Judah and S.
Shelah,} {\it The Cicho\'n diagram,} Journal of Symbolic
Logic, vol. 58 (1993), pp. 401-423.
\smallskip
\itemitem{[Fr]} {\capit D. Fremlin,} {\it Cicho\'n's diagram,}
S\'eminaire Initiation \`a l'Analyse (G. Choquet,
M. Rogalski, J. Saint Raymond), Publications Math\'ematiques
de l'Universit\'e Pierre et Marie Curie, Paris, 1984,
pp. 5-01 - 5-13. 
\smallskip
\itemitem{[Je]} {\capit T. Jech,} {\it Set theory,} Academic Press,
San Diego, 1978.
\smallskip
\itemitem{[Ku 1]} {\capit K. Kunen,} {\it Set theory,} North-Holland,
Amsterdam, 1980.
\smallskip
\itemitem{[Ku 2]} {\capit K. Kunen,} {\it Random and Cohen reals,}
Handbook of set--theoretic topology, K. Kunen and J. Vaughan
(editors), North--Holland, Amsterdam, 1984,
pp. 887-911
\smallskip
\itemitem{[Mi]} {\capit A. Miller,} {\it Arnie Miller's problem list,}
in: Set Theory of the Reals (H. Judah, ed.), Israel Mathematical
Conference Proceedings, vol. 6, 1993, pp. 645-654.
\smallskip

\vfill\eject\end